\documentclass[12pt,a4paper]{article}

\usepackage{latexsym}

\begin{document}

\renewcommand{\thefootnote}{\fnsymbol{footnote}}{}

\begin{center}

On the Jacobian conjecture

\footnote[0]{ {\it 2000 Mathmatics Subject Classification.} \quad 
Primary 14R15; Secondary 32M17, 32W50.
} 

\end{center}

\begin{center}

By

Y{\scriptsize UKINOBU} ADACHI

\end{center}

\vspace{5mm}

{\bf Abstract} \  \ We show that the Jacobian conjecture 
of the two dimensional case is true.

\vspace{5mm}

{\bf 0. Introduction}

\vspace{5mm}

Let $T = (P(x,y),Q(x,y)) : \mathbf{C}^2(x,y) \rightarrow \mathbf{C}^2(\xi,\eta)$ be a polynomial mapping 
with the Jacobian $J(P,Q) = 1$. Then we will show that $T \in Aut(\mathbf{C}^2)$ where $Aut(\mathbf{C}^2)$ is the
group of automorphisms of $\mathbf{C}^2$ (Theorem 3.3).
By virtue of S. Kaliman [5], in order to prove that $T \in Aut(\mathbf{C}^2)$, it suffices to prove this 
when the level curve $\{P = \alpha\}$ is irreducible for every $\alpha \in \mathbf{C}$.

From Lemma 2.3 [1], if the following linear partial differential equation where $u(x,y)$ is an unknown function

\[\frac{\partial(P,u)}{\partial(x,y)} = g(x,y) ,\        \ (1)\]

\noindent
has an entire solution for every entire function $g(x,y)$ of $\mathbf{C}^2$ with $P(x,y)$ being fixed as above,
the level curve $\{P = \alpha\}$ is simply connected for every $\alpha \in \mathbf{C}$.
It is now easy to see that $T \in Aut(\mathbf{C}^2)$ (cf. Theorem 3.2 in [1]).

It is to be noted that there is a polynomial $R(x,y)$ whose every level curve is nonsingular, irreducible and
not simply-connected by Bartolo, Cassou-Nogu\`es and Velasco [2].  Such an examle is shown in [4, p.256]. 
Thus, presupposing the existence of a polynomial $Q(x,y)$ with $J(P,Q) = 1$ is indispensable 
in solving the above key problem related to the equation (1) in Theorem 2.4.

In the following sections, it is assumed that polynomials $P(x,y)$ and $Q(x,y)$ satisfy 
the condition $J(P,Q) = 1$ and the level curve $\{P = \alpha\}$ is irreducible for 
every $\alpha \in \mathbf{C}$.  This assumption will be hereafter called the condition (K) also and expressed
in such a manner that $(P,Q)$ satisfies the condition (K), for example.

\vspace{5mm}

{\bf 1. Preliminary}

\vspace{5mm}

It is easy to see that the polynomial $Q(x,y)$ is  primitive, namely all of level curves are irreducible
except for a finite number of them because $Q(x,y)$ satisfies $J(P,Q) = 1$. It is obvious that $P(x,y)$ is also 
primitive. It is obvious also that every level curve of $Q(x,y)$ is nonsingular.

The following proposition is a summarized statement from Suzuki [7, p.242].

\vspace{5mm}

Proposition 1.1. \  \ {\it Except for a finite number of exceptional level curves $\{P = \alpha_i\}_{i= 1, \dots, k}$,
every general level curve $\{P = \alpha\}$ is the same type as $(g,n)$.
That is, the normalization of $\{P = \alpha\}$ has a genus $g$ and n-punctured boundaries.
If we take a sufficiently small disk $U(\alpha)$ with the center $\alpha$ which does not
contain $\{\alpha_1, \dots, \alpha_k\}$, the set $E_{\alpha} := \{(x,y); U(\alpha) \ni P(x,y)\}$
is considered to be a topologically trivial fiber space on $U(\alpha)$ of projection $P$
with a fiber $R_0$, where $R_0$ is a finite Riemann surface of $(g,n)$ type .}

\vspace{5mm}

Since $\{P = \alpha_i\}_{i=1, \dots, k}$ is a nonsingular curve, the following proposition, which is well known,
is a restatement of Lemma 1 in Nishino-Suzuki [6], for example.  

\vspace{5mm}

Proposition 1.2. \  \ {\it For every relatively compact and connected open set $K$ on $\{P = \alpha_i\}$, 
there is a relatively compact domain 
$M(K)$ in $\mathbf{C}^2$ such that $M(K) \subset E_{\alpha_i} := \{(x,y); U(\alpha_i) \ni
P(x,y)\}$ where $U(\alpha_i)$ is a sufficiently small disk with the center $\alpha_i$ and
$M(K)$ is considered to be a topologically trivial fiber space on $U(\alpha_i)$ of projection $P$
with a fiber $K$.}

\vspace{5mm}

Proposition 1.3. \  \ {\it For every $\alpha \in \mathbf{C}$, we consider 
a set $E_{\alpha} := \{(x,y); U(\alpha)
\ni P(x,y)\}$ as a fiber space such as $P(x,y) : E_{\alpha} \rightarrow 
U(\alpha)$. If we take a sufficiently small disk $U(\alpha)$, there is a 
global holomorphic section $L$, such as $\{Q(x,y) = \gamma\} \cap  E_{\alpha}$
where $\{Q(x,y)= \gamma\}$ is irreducible in $\mathbf{C}^2$ and $\{Q(x,y) = \gamma\}
\cap \{P(x,y) = \delta\} \neq \emptyset $ for every $\delta \in \mathbf{C}$.}

\vspace{5mm}

Proof. \    \ Let a point $(x_0,y_0)$ satisfy $P(x_0,y_0) = \alpha$ and $Q(x_0,y_0) = \beta$. Since $(P,Q)$ is locally
a biholomorphic map of near $(x_0,y_0)$ to near $(\alpha,\beta)$, $\{Q(x,y) = \beta'\}$ with $|\beta - \beta'| < \varepsilon$
, where $\varepsilon$ is sufficiently small positive number, is transversally intersects with $E_{\alpha}$ if we take a
sufficiently small disk $U(\alpha)$. Since $Q(x,y)$ is primitive as mentioning above, there is a complex number $\gamma$
such that $|\beta - \gamma| < \varepsilon$ and $\{Q(x,y) = \gamma\}$ is irreducible in $\mathbf{C}^2$.

It is well known that the image of $T$ is a Zariski open set. If the complement of the image of $T$ contains an irreducible algebraic
curve such as $\{C(x,y) = 0\}$, then $C(P,Q)$ is a polynomial and $C(P,Q) \neq 0$. This contradicts to the fact that $T$ is
a nondegenerate map. Therefore the complement of the image of $T$ is a set of finite number of points of $\mathbf{C}^2(\xi,\eta)$.
And for almost $\gamma$ with $|\beta - \gamma| < \varepsilon$, $\{Q(x,y) = \gamma\} \cap \{P(x,y) = \delta\} \neq \emptyset$
for every value $\delta \in \mathbf{C}$. Since $Q(x,y)$ is a primitive polynomial, there is a $\gamma$ which satisfies a 
condition of the proposition. \    \ $\Box$

\vspace{5mm}

Following proposition is easy to see.

\vspace{5mm}

Proposition 1.4. \  \ {\it Let $\alpha \in \mathbf{C} - \{\alpha_1, \dots ,\alpha_k\}$. Let $U(\alpha)$
and $E_{\alpha}$ be defined in the same way as Proposition 1.1 and have a global holomorphic section $L$ 
in the same way as  Proposition 1.3.
Then the universal covering space $\tilde{E_{\alpha}}$ of $E_{\alpha}$ whose base points are points of $L$
can be thought  as $U(\alpha) \times \tilde{R_0}$ where $\tilde{R_0}$ is a universal covering surface of $R_0$.
Roughly speaking, $\tilde{E_{\alpha}}$ is a fiber space whose fibers are universal covering surfaces for every fiber 
of $E_{\alpha}$.}

\vspace{5mm}

{\bf 2. Construction of a global solution of the differential
equation (1)}

\vspace{5mm}

We consider a linear partial differential equation

\[\frac{\partial(P,u)}{\partial(x,y)} = g(x,y), \        \ (1)\]

\noindent
where $u$ is an unknown function and $g(x,y)$ is an arbitrary entire function of $\mathbf{C}^2$.
We assume that there is a polynomial $Q(x,y)$ such that $(P,Q)$  satisfies the condition (K).

\vspace{5mm}

Since (1) is a linear and nonsingular partial differential equation, 
there is a unique local single-valued holomorphic solution $u$ 
near $L$, where $L$ is the same as in Proposition 1.3,
having any given holomorphic initial data on $L$ by Cauchy-Kowalevskaya's theorem (If necessary,
we can take $U(\alpha)$ to be smaller). For applying Cauchy-Kowalevskaya's theorem we remark that $\{P = \alpha\}$ is the characteristic curve of (1)
and a nonsingular curve $L$ is transversely crossing to $\{P = \alpha\}$.

\vspace{5mm}

Lemma 2.1. \  \ {\it Let $P(x,y) : E_{\alpha} \rightarrow U(\alpha)$ be the same as that of Proposition 1.4 and
$u(x,y)$ be the single-valued holomorphic solution of $(1)$ near $L$ having arbitrarily given holomorphic initial data on $L$.
Then, $u(x,y)$ has an analytic continuation along any path in $E_{\alpha}$ with any starting point on $L$.}

\vspace{5mm}

Proof. \  \ Since the characteristic curve of the equation (1) is $\{P = \alpha\}$ and the equation (1) can be
thought  as an analytic family of holomorphic 1-forms such as $du = \frac{g}{P_y} dx = - \frac{g}{P_x} dy$
where $P_x dx + P_y dy = 0$ (see p.637 in [1]), $u(x,y)$ is uniquely fixed as a multi-valued
integral of a holomorphic 1-form for every fiber of $E_{\alpha}$.

From Proposition 1.4, $u(x,y)$ is determined on $\tilde{E_{\alpha}}$ as a single-valued function.
As $u(x,y)$ is a holomorphic solution of two variables near $L$ and holomorphic on every fiber as 
a function of one complex variable,
$u(x,y)$ is a single-valued holomorphic solution on $\tilde{E_{\alpha}}$ by the well-known Hartogs theorem.

Then, $u(x,y)$ has an analytic continuation along any path in $\tilde{E_{\alpha}}$ with the starting point of $L$.
As any path in $E_{\alpha}$ with the starting point of $L$ can be considered as the one in $\tilde{E_{\alpha}}$,
the above lemma holds true. \    \ $\Box$

\vspace{5mm}

Lemma 2.2. \  \ {\it For $\{P = \alpha_i\}, i \in \{1, \dots ,k\}$, we take a disk $U(\alpha_i)$ sufficiently
small such that $U(\alpha_i) \cap \{\alpha_1, \dots, \alpha_{i-1}, \alpha_{i+1}, \cdots ,\alpha_k \} = \emptyset$,
$E_{\alpha_i} := \{(x,y); U(\alpha_i) \ni P(x,y)\}$ has a global holomorphic section $L$ as in Proposition 1.3 and there exists a
holomorphic solution $u(x,y)$ near $L$ having arbitrarily given holomorphic initial data on $L$.
Let $U^*(\alpha_i) = U(\alpha_i) - \{\alpha_i\}, L^* = L - \{P = \alpha_i\}$ and $E_{\alpha_i}^* = \{(x,y); U^*(\alpha_i)
\ni P(x,y)\}$. 

Then, $u(x,y)$ can be analytically continued along any path $l_{rs}$ which starts from a point $r$ of $L^*$ and extends to any point $s$
such that $l_{rs} \subset E_{\alpha_i}^*$.}

\vspace{5mm}

Proof. \  \ Let  $\delta_j$ be a disk in $U^*(\alpha_i)$ such that 
$\bigcup_{j=1,2, \cdots} \delta_j = U^*(\alpha_i)$ and $E_{\delta_j} := \{(x,y); \delta_j \ni P(x,y)\}$ is regarded as a 
topologically trivial fiber space. 

Let $r \in E_{\delta_j} \cap L^*$. We prove at  first for the case $s \in E_{\delta_k}$, $\delta_j \cap \delta_k \neq
\emptyset $ and $l_{rs} \subset E_{\delta_j} \cup E_{\delta_k}$.

By the (topological) triviality of $E_{\delta_j}$ and $E_{\delta_k}$, the curve $l_{rs}$ is homotopic to $l_{rr'} \cup l_{r's'}
\cup l_{s's}$ in $E_{\delta_j} \cup E_{\delta_k}$ where $r' \in E_{\delta_j} \cap E_{\delta_k} \cap L^*$, $l_{rr'} \subset  L^*$, $l_{r's'} \subset
\{P = P(r')\}$ and $l_{s's} \subset E_{\delta_k}$.

From Lemma 2.1, $u(x,y)$, which has given initial data near $r' \in L^*$, can be analytically continued along $l_{r's'} \cup
l_{s's}$ and  $u(x,y)$ has a function element $u_1$ near $s$. As $l_{rs} \simeq l_{rr'} \cup l_{r's'} \cup l_{s's}$ and
$u(x,y)$ has an analytic continuation along any path in $E_{\delta_j}$ and $E_{\delta_k}$ respectively by Lemma 2.1,
$u(x,y)$ can be analytically continued along $l_{rs}$ and has a function element $u_1$ near $s$.

Next we prove for the case $s \in E_{\delta_l}, \delta_j \cap \delta_k \neq \emptyset, \delta_k \cap \delta_l
\neq \emptyset$ and $l_{rs} \subset E_{\delta_j} \cup E_{\delta_k} \cup E_{\delta_l}$.
To simplyfy the proof (without loss of generality), we assume that $l_{rs} = l_{rs'} \cup l_{s's''} \cup l_{s''s}$
where $l_{rs'} \subset E_{\delta_j} \cup E_{\delta_k}, l_{s's''} \subset E_{\delta_k} \cup E_{\delta_l}$ and
$l_{s''s} \subset E_{\delta_l}$ by the triviality of $E_{\delta_j}, E_{\delta_k}$ and $E_{\delta_l}$.

By the triviality of $E_{\delta_j}$ and $E_{\delta_k}$, $l_{rs'} \simeq l_{rr'} \cup l_{r's'}$ where 
$r' \in L^* \cap E_{\delta_j} \cap E_{\delta_k}, l_{rr'} \subset L^*$ and $l_{r's'} \subset \{P = P(r') = P(s')\}$. 
By the triviality of $E_{\delta_k}$ and $E_{\delta_l},
l_{rs'} \cup l_{s's''} \simeq l_{rr'} \cup l_{r's'} \cup l_{s's''} \simeq l_{rr'} \cup l_{r'r''} \cup l_{r''s''}$
where $r'' \in L^* \cap E_{\delta_k} \cap E_{\delta_l}$, $l_{r'r''} \subset L^*$ 
and $l_{r''s''} \subset \{P = P(r'') = P(s'')\}$. 
Then $l_{rs} \simeq l_{rr'} \cup l_{r'r''} \cup l_{r''s''} \cup l_{s''s}$ where $l_{s''s} \subset E_{\delta_l}$.
For the same reason as above, $u(x,y)$ can be analytically continued along $l_{rs}$.
Based on the succesive discussions as above, the above lemma is proved easily. \    \ $\Box$

\vspace{5mm}

Lemma 2.3. \  \ {\it Let $E_{\alpha_i}$ be the same of Lemma 2.2. Then, $u(x,y)$ can be analytically continued along
any path $l_{rs}$ which starts a point $r$ of $L$ to every point $s$ of $E_{\alpha_i}$ where $l_{rs} \subset E_{\alpha_i}$.}

\vspace{5mm}

Proof. \  \ From the remark made at the outset of this section and with the same reasoning used for
the proof of Lemma 2.1, we remark that $u(x,y)$ can be analytically continued along any path on $\{P = \alpha_i\}$ 
as a multi-valued holomorphic function of two variables having given initial data near $L$.

Let $r \in L, s \in E_{\alpha_i}$ and $l_{rs} = ((x(t),y(t)), t \in [0,1]$, where $r = (x(0),y(0))$ and $s = (x(1),y(1))$. Then $\eta(t) = Q(x(t),y(t))$ is
a continuous path on $\eta$-plane. We lift up $\eta(0)$ to $r_0 \in \{P = \alpha_i\} \cap L$ by $T^{-1}$ 
and lift up continuously
$\eta(t)$ to the point of $\{P = \alpha_i\}$ by $T^{-1}$. We set such a path to $l_{r_0 s_0}$.
More precisely speaking, $\eta = Q(x,y)_{|\{P=\alpha_i\}}$ takes every value at same times except at most finite values $\gamma_1, \gamma_2, \cdots , \gamma_m$.
The map $T^{-1}$ is able to continue holomorphically along any path from $r_0$ on $\{P = \alpha_i\}$ avoiding finite points.

From Proposition 1.2, if we take a relatively compact connected open set $K$ on $\{P = \alpha_i\}$ 
in which $l_{r_0s_0}$ is contained as a relatively compact subset,
we can take a tubular neighborhood $M(K)$ of $K$. 
Then, the solution $u(x,y)$ of (1) having given initial data on $L \cap M(K)$ can be analytically continued
along every path in $M(K)$ if we take $M(K)$ sufficiently thin with the same reasoning as in Lemma 2.1 using Lemma 1.2.

There are two cases to be considered.

(1) If $s \in E_{\alpha_i}^*$, the path $l_{rs}$ is homotopic to $l_{rr'} \cup l_{r's'} \cup l_{s's}$ where
$r' \in L^* \cap M(K), l_{rr'} \subset L , l_{r's'} \subset M(K) - \{P = \alpha_i\}$ and $l_{s's} 
\subset E_{\alpha_i}^*$ in the same way as one in the proof of Lemma 2.2. From Lemma 2.2, $u(x,y)$ can be analytically
continued along $l_{rr'} \cup l_{r's'} \cup l_{s's}$ and $u(x,y)$ has a function element $u_1$ near $s$. 
Then $u(x,y)$ can be analytically continued along $l_{rs}$ and $u(x,y)$ has a function element $u_1$ near $s$.

(2) If $s \in \{P = \alpha_i\}$, the path $l_{rs}$ is homotopic to $l_{rr'} \cup l_{r's}$ where
$r' \in L \cap \{P = \alpha_i\}, l_{rr'} \subset L, l_{r's} \subset \{P = \alpha_i\} \cap M(K)$.
From the property of $M(K)$, $u(x,y)$ can be analytically continued along $l_{rr'} \cup l_{r's}$ and $u(x,y)$ has a function
element $u_1$ near $s$. Then $u(x,y)$ can be analytically continued along $l_{rs}$ and $u(x,y)$ has a function element $u_1$
near $s$. \    \ $\Box$

\vspace{5mm}

Theorem 2.4. \  \ {\it For every entire function $g(x,y)$ of $\mathbf{C}^2$, the equation $(1)$ has an entire solution $u(x,y)$.}

\vspace{5mm}

Proof. \  \ Suppose that $P(x,y) : E_{\alpha} \rightarrow U(\alpha)$ satisfies Lemma 2.1 and suppose that $P(x,y) : E_{\alpha_i} \rightarrow 
U(\alpha_i)$ satisfies Lemma 2.3. We take $\{U(\alpha)\}$ and $\{U(\alpha_i)\}$ smaller than themselves, if necessary, 
to satisfy Lemma 2.2 and vary notations of them altogether to be a countable covering $\{U_j\}_{j=1,2,\cdots}$ of $\mathbf{C}$. 
We denote $E_j$ by $\{(x,y); U_j \ni P(x,y)\}$.
Let $L_j$ be a global holomorphic section as in Proposition 1.3 of $P(x,y) : E_j
\rightarrow U_j$ and let $u_j$ be a solution of Eqaution (1) having given an initial holomorphic data on $L_j$ which can be continued analytically
along every path in $E_j$. Let $u_j^0$ is a branch of $u_j$ which coincides to  given an initial data on $L_j$.

Now we proceed to prove Theorem 2.4 in the following four steps (1) to (4).

(1) At first, consider a case where $U_j \cap U_k \neq \emptyset$.
We take $r_j \in L_j \cap E_j \cap E_k$ and $r_k \in L_k \cap E_k \cap E_j$ and fix a path $l_{jk}$ from $r_j$ to
$r_k$ such that $l_{jk} \subset E_j \cap E_k$ and $l_{kj} = - l_{jk}$.
We have  $u_j^0$  continue analytically along $l_{jk}$ from $r_j$ to $r_k$ and have it  continue  analytically on every point of 
$L_k \cap E_j$ and we denote such a single-valued branch by $u_j^k$ since $L_k \cap E_j$ is simply connected.
As $J(P,u_k^0 - u_j^k) = 0$ in the neighborhood of $L_k \cap E_j$, $u_k^0 - u_j^k = \varphi_{jk}(P(x,y))$
where $\varphi_{jk}$ is a single-valued holomorphic function on $L_k \cap E_j$ and then it is a single-valued holomorphic 
one on $U_j \cap U_k$.
Let $u_k^0$ continue analytically along $l_{kj}$ from $r_k$ to $r_j$ and let it continue analytically on every point of $L_j \cap
E_k$ and we denote such a single-valued branch by $u_k^j$ since $L_j \cap E_k$ is simply connected.
Similar to the above, $u_j^0 - u_k^j = \varphi_{kj}(P(x,y))$, where $\varphi_{kj}$ can be considered as a single-valued 
holomorphic function on $U_k \cap U_j$.
If we have $u_j^0 - u_k^j$ continue analytically along $l_{jk}$ from $r_j$ to $r_k$ and have it continue analytically on $L_k \cap E_j$,
then its branch on $L_k \cap E_j$ is $u_j^k - u_k^0$. Hence, $\varphi_{jk} = - \varphi_{kj}$ on $U_j \cap U_k$.

(2) When $U_j \cap U_k \cap U_l \neq \emptyset$, we will prove that $\varphi_{jk} + \varphi_{kl} + \varphi_{lj} = 0$
on $U_j \cap U_k \cap U_l$ where $\varphi_{jk}, \varphi_{kl}$ and $\varphi_{lj}$ are determined in the same way as with the
case (1). Since $\varphi_{jk}(P) = u_k^0 - u_j^k = u_k^j - u_j^0$ on $L_j \cap E_k$, $\varphi_{kl}(P) = u_l^0 - u_k^l
= u_l^j - u_k^j$ on $L_j \cap E_k \cap E_l$ and $\varphi_{lj}(P) = u_j^0 - u_l^j$ on $L_j \cap E_l$, $\varphi_{jk} + \varphi_{kl} +
\varphi_{lj} = 0$ on $U_j \cap U_k \cap U_l$.

Since $\{\varphi_{jk} ; U_j \cap U_k \neq \emptyset\}$ is a Cousin I data from the fact of (1) and above, and
$\mathbf{C}$ is a Cousin I domain, there exist $\varphi_j \in \mathcal{O}(U_j)$ and $\varphi_k \in \mathcal{O}(U_k)$
such that $\varphi_{jk} = \varphi_j - \varphi_k$ on $U_j \cap U_k \neq \emptyset$.

(3) We will prove that $\{u_j^0 + \varphi_j(P(x,y))$ on $E_j\}_{j=1,2,\cdots}$ is a connected multi-valued function $u(x,y)$.
When $E_j \cap E_k \neq \emptyset$, we have $u_j^0 + \varphi_j(P(x,y))$ continue analytically along $l_{jk}$ to $L_k \cap
E_j$. Then, near $L_k \cap E_j$

$u_k^0 + \varphi_k(P(x,y)) - (u_j^k + \varphi_j(P(x,y))$

$= u_k^0 - u_j^k - \{\varphi_j(P(x,y)) - \varphi_k(P(x,y))\}$

$= \varphi_{jk}(P(x,y)) - \varphi_{jk}(P(x,y)) = 0$.

\noindent By successively applying the same reasoning as above, we prove the above fact.

(4)Lastly, we will prove that such multi-valued function $u(x,y)$ can be analytically continued along every path $l_{rs}$ in
$\mathbf{C}^2$ where $r \in L_j$ and $s \in E_k$ where $j$ and $k$ are arbitrary positive integers. 
We consider a new $u(x,y)$ where all function elements of $\{u_j^0 + \varphi_j(P(x,y)$ on $E_j\}_{j=1,2,\cdots}$ 
are continued analytically to every path in $\mathbf{C}^2$ if possible.

It is easy to see that $l_{rs} = l_{rs_1} \cup l_{s_1s_2} \cup \cdots \cup l_{s_ns}$ where
$l_{rs_1} \subset E_j \cup E_{k_1}, l_{s_1s_2} \subset E_{k_1} \cup E_{k_2}, s_1 \in E_{k_1},
\cdots, l_{s_ns} \subset E_{k_n} \cup E_k, s_n \in E_{k_n}$.
By the same way of the proof of Lemma 2.3 using $Q(x,y)$ and $T^{-1}$, $l_{rs_1}$ is homotopic to
$l_{rr'} \cup l_{r's'} \cup l_{s's''} \cup l_{s''s_1}$, where $l_{rr'} \in L_j,
l_{r's'} \subset \{P(x,y) = P(r')\} \subset E_j \cap E_{k_1}, s' \in L_{k_1},
l_{s's''} \subset L_{k_1}$ and $l_{s''s_1} \subset E_{k_1}$.
Since $u(x,y)$ can be continued analytically along $l_{rr'} \cup l_{r's'} \cup l_{s's''} \cup l_{s''s_1}$
by Lemma 2.1 and Lemma 2.3. By sucessively applying the same reasoning as above, we can continue $u(x,y)$ analytically
along $l_{rs}$.

Therefore, such $u(x,y)$ is a connected multi-valued function 
and can be continued analytically along every path in $\mathbf{C}^2$. From the monodromy theorem, $u(x,y)$ is a 
single-valued entire function. \    \ $\Box$

\vspace{5mm}

{\bf 3. Conclusion}

\vspace{5mm}

Lemma 3.1 (cf. Lemma 2.3 in [1]). \  \ {\it The level curve $\{P = \alpha\}$ for every $\alpha \in \mathbf{C}$
is simply-connected.}

\vspace{5mm}

Proof. \    \ If $S_0 := \{P = \alpha_0\}$ is not simply-connected, there is a holomorphic 1-form $a(x)dx$
or $b(y)dy$, where $b(y) = a(x)(- \frac{P_y}{P_x})$ whose integral on $S_0$ is a multi-valued function by Behnke-Stein
theorem in [3]. If we set $a(x)P_y = - b(y)P_x$, it represents a holomorphic function $g_0$ on $S_0$.
As a consequence of Cartan's Theorem B, there is an entire function $g$ of $\mathbf{C}^2$ such that $g|_{S_0} = g_0$.
It is easy to see that the equation (1) for such $g$ never has any single-valued entire solution.
This contradicts Theorem 2.4. \    \ $\Box$

\vspace{5mm}

Theorem 3.2 (cf. Theorem 3.2 in [1]). \  \ {\it If the polynomial map $T = (P(x,y),Q(x,y))$ satisfies the condition $(K)$,
then $T \in Aut(\mathbf{C}^2)$.}

\vspace{5mm}

Proof. \  \ If we restrict $Q(x,y)$ to $S := \{P(x,y) = \alpha\}$, $dQ$ has a pole of order 2 at $\infty$ because the degree of $dQ =
-2$ by the Riemann-Roch theorem and does not vanish on $S$. Then $P|_S$ takes every value once by the residue theorem.
It is true for every $\alpha \in \mathbf{C}$ and therefore $T \in Aut(\mathbf{C}^2)$. \     \ $\Box$

\vspace{5mm}

By virtue of Kaliman [5], the Jacobian conjecture of the two dimensional case is proved, that is,

\vspace{5mm}

Theorem 3.3. \  \ {\it If the polynomial map $T = (P(x,y),Q(x,y))$ satisfies $J(P,Q) = 1$, then $T \in Aut(\mathbf{C}^2)$.}

\vspace{5mm}

\begin{center}

{\bf References}

\end{center}

\vspace{5mm}

[1] Y. Adachi, Condition for global exsistence of holomorphic solutions of a certain differential equation on 
a Stein domain of $\mathbf{C}^{n+1}$ and its applications, J. Math. Soc. Japan, {\bf 53} (2001), 633-644.

[2] E.A. Bartolo, P. Cassou-Noug\`es and I.L. Velasco, On polynomials whose fibers are irreducible with no critical points, 
Math. Ann., {\bf 299} (1994), 477-490.

[3] H. Behnke and K. Stein, Entwicklung analytischen Functionen auf Riemannshen Fl\"achen, Math. Ann., {\bf 120} (1949),
430-461.

[4] A. van den Essen, Polynomial Automorphisms and the Jacobian Conjecture, Birkh\"auser,2000.

[5] S. Kaliman, On the Jacobian conjecture, Proc. Amer. Math. Soc., {\bf 117} (1993), 45-51.

[6] T. Nishino and M. Suzuki, Sur les singularit\'es essentielles et isol\'ees des applications holomorphes \`a valeurs
dans une surface complexe, Publ. Res. Inst. Math. Sci. Kyoto Univ., {\bf 16} (1980), 461-497.

[7] M. Suzuki, Propri$\acute{e}$t$\acute{e}$s topologiques des polyn$\hat{o}$mes de deux variables complexes, et automorphismes
alg$\acute{e}$briques de l'espace $\mathbf{C}^2$, J. Math. Soc. Japan, {\bf 26} (1974), 241-257.

\vspace{10mm}

Yukinobu Adachi

2-12-29, Kurakuen,

Nishinomiya City, Hyogo 662-0082, Japan

E-mail: fwjh5864@nifty.com

\end{document}